\documentclass[11pt,twoside,a4paper]{article}

\usepackage{setspace}

\usepackage{amssymb}
\usepackage{amsmath}
\usepackage[dvips]{graphicx}
\usepackage[all,2cell,ps]{xy}
\title{Higher Minors and Van Kampen's Obstruction}
\author{Eran Nevo \footnote{Institute of Mathematics, Hebrew
University, Jerusalem Israel, E-mail address:
eranevo@math.huji.ac.il}}
 \newtheorem{thm}{Theorem}[section]
 \newtheorem{cor}[thm]{Corollary}
 \newtheorem{lem}[thm]{Lemma}
\newtheorem{defin}[thm]{Definition}
 \newtheorem{prop}[thm]{Proposition}
\newtheorem{prob}[thm]{Problem}
\newtheorem{conj}[thm]{Conjecture}
\newtheorem{ex}[thm]{Example}
\newtheorem{obs}[thm]{Observation}

\begin{document}
\maketitle
%\doublespacing

\begin{abstract}
We generalize the notion of graph minors to all (finite)
simplicial complexes. For every two simplicial complexes $H$ and
$K$ and every nonnegative integer $m$, we prove that if $H$ is a
minor of $K$ then the non vanishing of Van Kampen's obstruction in
dimension $m$ (a characteristic class indicating non embeddability
in the $(m-1)$-sphere) for $H$ implies its non vanishing for $K$.
As a corollary, based on results by Van Kampen \cite{VanKampen}
and Flores \cite{Flores}, if $K$ has the $d$-skeleton of the
$(2d+2)$-simplex as a minor, then $K$ is not embeddable in the
$2d$-sphere.

We answer affirmatively a problem asked by Dey et. al. \cite{Dey}
concerning topology-preserving edge contractions, and conclude
from it the validity of the generalized lower bound inequalities
for a special class of triangulated spheres.
\end{abstract}

\section{Introduction}\label{secIntro}
The concept of graph minors has proved be to very fruitful.
% in particular it is used to detect embeddability of graphs in
%surfaces.
A famous result by Kuratowski asserts that a graph can be embedded
into a $2$-sphere if and only if it contains neither of the graphs
$K_5$ and $K_{3,3}$ as minors. We wish to generalize the notion of
graph minors to all (finite) simplicial complexes in a way that
would produce analogous statements for embeddability of higher
dimensional complexes in higher dimensional spheres. We hope that
these higher minors will be of interest in future research, and
indicate some results and problems to support this hope.
\newline
Let $K$ and $K'$ be simplicial complexes. $K\mapsto K'$ is called a
\emph{deletion} if $K'$ is a subcomplex of $K$. $K\mapsto K'$ is
called an \emph{admissible contraction} if $K'$ is obtained from $K$
by identifying two distinct vertices of $K$, $v$ and $u$, such that
$v$ and $u$ are not contained in any missing face of $K$ of
dimension $\leq \rm{dim}(K)$. (A set $T$ is called a missing face of
$K$ if it is not an element of $K$ while all its proper subsets
are.) Specifically, $K'=\{T: u\notin T \in K\}\cup\{(T\setminus
\{u\})\cup \{v\}: u\in T \in K\}$. An equivalent formulation of the
condition for admissible contractions is that the following holds:
\begin{equation}\label{eqSkelLinkCond}
\rm{skel}_{\rm{dim}(K)-2}(\rm{lk}(v,K) \cap \rm{lk}(u,K)) =
\rm{lk}(\{v,u\},K)
\end{equation}
where $\rm{skel}_{m}(K)$ is the subcomplex of $K$ consisting of
faces of dimension $\leq m$ and $\rm{lk}(v,K):=\{F\in K: v\notin F,
F\cup \{v\}\in K\}$ is the link of $v$ in $K$. For $K$ a graph,
(\ref{eqSkelLinkCond}) just means that $\{v,u\}$ is an edge in $K$.

We say that a simplicial complex $H$ is a \emph{minor} of $K$, and
denote it by $H<K$, if $H$ can be obtained from $K$ by a sequence of
admissible contractions and deletions (the relation $<$ is a partial
order). Note that for graphs this is the usual notion of a minor.
\newline
\textbf{Remarks}: (1) In equation (\ref{eqSkelLinkCond}), the
restriction to the skeleton of dimension at most $\rm{dim}(K)-2$ can
be replaced by restriction to the skeleton of dimension at most
$\rm{min}\{\rm{dim}(\rm{lk}(u,K)),\rm{dim}(\rm{lk}(v,K))\}-1$,
making the condition for admissible contraction \emph{local}, and
weaker. All the results and proofs in this paper hold verbatim for
this notion of a minor as well.

(2) In the definition of a minor, without loss of generality we may
replace the local condition from the remark above by the following
stronger local condition, called the \emph{Link Condition} for
$\{u,v\}$:
\begin{equation}\label{eqLinkCond}
\rm{lk}(u,K)\cap \rm{lk}(v,K)=\rm{lk}(\{u,v\},K).
\end{equation}
To see this, let $K\mapsto K'$ be an admissible contraction which is
obtained by identifying the vertices $u$ and $v$ where
$\rm{dim}(\rm{lk}(u,K))\leq \rm{dim}(\rm{lk}(v,K))$. Delete from $K$
all the faces $F\uplus\{u\}$ such that $F\uplus\{u,v\}$ is a missing
face of dimension $\rm{dim}(\rm{lk}(u,K))+2$, to obtain a simplicial
complex $L$. Note that $\{u,v\}$ satisfies the Link Condition in
$L$, and the identification of $u$ with $v$ in $L$ results in $K'$.
I thank an anonymous referee for this remark.

We now relate this minor notion to Van Kampen's obstruction in
cohomology; following Sarkaria \cite{SarkariaMax} we will work with
deleted joins and with $\mathbb{Z}_2$ coefficients (background and
definitions appear in the next section).
\begin{thm}\label{thm o(L)}
Let $\rm{Sm}^m(L)\in H^{m}_S(L_*,\mathbb{Z}_2)$ denote Van
Kampen's obstruction (in equivariant cohomology) for a simplicial
complex $L$, where $L_*$ is the deleted join of $L$. Let $H$ and
$K$ be simplicial complexes. If $H<K$ and $\rm{Sm}^m(H)\neq 0$
then $\rm{Sm}^m(K)\neq 0$.
\end{thm}
For any positive integer $d$ let $H(d)$ be the $(d-1)$-skeleton of
the $2d$-dimensional simplex. A well known result by Van Kampen and
Flores \cite{Flores,VanKampen} asserts that the Van Kampen
obstruction of $H(d)$ in dimension $(2d-1)$ does not vanish, and
hence $H(d)$ is not embeddable in the $2(d-1)$-sphere (note that the
case $H(2)=K_5$ is part of the easier direction of Kuratowski's
theorem).
\begin{cor}\label{cor Emb}
For every $d\geq 1$, if $H(d)<K$ then $K$ is not embeddable in the
$2(d-1)$-sphere. $\square$
\end{cor}
\textbf{Remark}:\label{rem Emb} Corollary \ref{cor Emb} would also
follow from the following conjecture:
\begin{conj}\label{conjMinorEmb}
If $H<K$ and $K$ is embeddable in the $m$-sphere then $H$ is
embeddable in the $m$-sphere.
\end{conj}

The following theorem answers in the affirmative a question asked by
Dey et. al. \cite{Dey}, who already proved the dimension $\leq 3$
case.
\begin{thm}\label{thmDey}
Given an edge in a triangulation of a compact PL (piecewise
linear)-manifold without boundary, its contraction results in a
PL-homeomorphic space if and only if it satisfies the Link Condition
(\ref{eqLinkCond}).
\end{thm}

In Section \ref{secBackground} we give the needed background on Van
Kampen's obstruction and Smith characteristic class. In Section
\ref{secProof} we prove Theorem \ref{thm o(L)} and show some
applications. In Section \ref{secOverZ} we prove an analogue of
Theorem \ref{thm o(L)} for deleted products over $\mathbb{Z}$. In
Section \ref{secTop} we prove Theorem \ref{thmDey} and deduce from
it some $f$-vector consequences. In Section \ref{secVersus} we
compare higher minors with graph minors.

\section{Algebraic-topological background}\label{secBackground}
The presentation here is based on work of Sarkaria
\cite{SarkariaMax,SarkariaUnpub} who attributes it to Wu \cite{Wu}
and all the way back to Van Kampen \cite{VanKampen}. It is a Smith
theoretic interpretation of Van Kampen's obstructions.

Let $K$ be a simplicial complex. The join $K*K$ is the simplicial
complex $\{S^1\uplus T^2: S,T\in K\}$ (the superscript indicates
two disjoint copies of $K$). The \emph{deleted join} $K_*$ is the
subcomplex $\{S^1\uplus T^2: S,T\in K, S\cap T=\emptyset\}$. The
restriction of the involution $\tau:K*K\longrightarrow K*K$,
$\tau(S^1\cup T^2)=T^1\cup S^2$ to $K_*$ is into $K_*$. It induces
a $\mathbb{Z}_2$-action on the cochain complex
$C^*(K_*;\mathbb{Z}_2)$. For a simplicial cochain complex $C$ over
$\mathbb{Z}_2$ with a $\mathbb{Z}_2$-action $\tau$, let $C_S$ be
its subcomplex of \emph{symmetric cochains}, $\{c\in C:
\tau(c)=c\}$. Restriction induces an action of $\tau$ as the
identity map on $C_S$. Note that the following sequence is exact
in dimensions $\geq 0$:
$$0\longrightarrow C_S(K_*)\longrightarrow C(K_*)\stackrel{id+\tau}{\longrightarrow}
C_S(K_*)\longrightarrow 0$$ where $C_S(K_*)\longrightarrow C(K_*)$
is the trivial injection. (The only part of this statement that may
be untrue for a non-free simplicial cochain complex $C$ over
$\mathbb{Z}_2$ with a $\mathbb{Z}_2$-action $\tau$, is that
${id+\tau}$ is surjective.) Thus, there is an induced long exact
sequence in cohomology
$$H_S^0(K_*)\stackrel{\rm{Sm}}{\longrightarrow} H_S^1(K_*)\longrightarrow...\longrightarrow
H_S^q(K_*)\longrightarrow H^q(K_*) \longrightarrow
H_S^q(K_*)\stackrel{\rm{Sm}}{\longrightarrow} H_S^{q+1}(K_*)
\longrightarrow... .$$ Composing the connecting homomorphism
$\rm{Sm}$ $m$ times we obtain a map $\rm{Sm}^m:
H_S^0(K_*)\longrightarrow H_S^m(K_*)$. For the fundamental
$0$-cocycle $1_{K_*}$, i.e. the one which maps $\sum_{v\in
(K_*)_0}a_vv\mapsto \sum_{v\in (K_*)_0}a_v \in \mathbb{Z}_2$, let
$[1_{K_*}]$ denotes its image in $H_S^0(K_*)$.
$\rm{Sm}^m([1_{K_*}])$ is called the $m$-th \emph{Smith
characteristic class} of $K_*$, denoted also as $\rm{Sm}^m(K)$.

\begin{thm}(Sarkaria \cite{SarkariaUnpub} Theorem 6.5, see also Wu \cite{Wu} pp.114-118.)\label{thmSmithH(d)}
For every $d\geq 1$, $\rm{Sm}^{2d-1}(1_{H(d)_*})\neq 0$.
\end{thm}

\begin{thm}(Sarkaria \cite{SarkariaUnpub} Theorem 6.4 and \cite{SarkariaMax} p.6)\label{thmSmith->nonEmb}
If a simplicial complex $K$ embeds in $\mathbb{R}^m$ (or in the
$m$-sphere) then $\rm{Sm}^{m+1}(1_{K_*})=0$.
\end{thm}
$Sketch\ of\ proof$: The definition of Smith class makes sense for
singular homology as well; the obvious map from the simplicial chain
complex to the singular one induces an isomorphism between the
corresponding Smith classes. The definition of deleted join makes
sense for subspaces of a Euclidean space as well (see e.g.
\cite{Matousek}, 5.5); thus an embedding $|K|$ of $K$ into
$\mathbb{R}^m$ induces a continuous $\mathbb{Z}_2$-map from $|K|_*$
into the join of $\mathbb{R}^m$ with itself minus the diagonal,
which is $\mathbb{Z}_2$-homotopic to the antipodal $m$-sphere,
$S^m$. The equivariant cohomology of $S^m$ over $\mathbb{Z}_2$ is
isomorphic to the ordinary cohomology of $\mathbb{R}P^m$ over
$\mathbb{Z}_2$, which vanishes in dimension $m+1$. We get that
$\rm{Sm}^{m+1}(S^m)$ maps to $\rm{Sm}^{m+1}(1_{|K|_*})$ and hence
the later equals to zero as well. But $|K_*|$ and $|K|_*$ are
$\mathbb{Z}_2$-homotopic, hence $\rm{Sm}^{m+1}(1_{K_*})=0$.
$\square$

\section{A proof of Theorem \ref{thm o(L)}}\label{secProof}
The idea is to define an injective chain map $\phi:
C_*(H;\mathbb{Z}_2)\longrightarrow C_*(K;\mathbb{Z}_2)$ which
induces $\phi(\rm{Sm}^m(1_{K_*}))=\rm{Sm}^m(1_{H_*})$ for every
$m\geq 0$.

\begin{lem}\label{lemAdContr->InjChainMap}
Let $K\mapsto K'$ be an admissible contraction. Then it induces an
injective chain map $\phi: C_*(K';\mathbb{Z}_2)\longrightarrow
C_*(K;\mathbb{Z}_2)$.
\end{lem}
$Proof$: Fix a labeling of the vertices of $K$, $v_0,v_1,..,v_n$,
such that $K'$ is obtained from $K$ by identifying $v_0\mapsto v_1$
where $\rm{dim}(\rm{lk}(v_0,K))\leq \rm{dim}(\rm{lk}(v_1,K))$.

Let $F\in K'$. If $F\in K$, define $\phi(F)=F$. If $F\notin K$,
define $\phi(F)=\sum\{(F\setminus v)\cup v_0: v\in F, (F\setminus
v)\cup v_0\in K\}$. Note that if $F\notin K$ then $v_1\in F$ and
$(F\setminus v_1)\cup v_0 \in K$, so the sum above is nonzero.
Extend linearly to obtain a map $\phi:
C_*(K';\mathbb{Z}_2)\longrightarrow C_*(K;\mathbb{Z}_2)$.

First, let us check that $\phi$ is a chain map, i.e. that it
commutes with the boundary maps $\partial$. It is enough to verify
this for the basis elements $F$ where $F\in K'$. If $F\in K$ then
$supp(\partial F)\subseteq K$, hence $\partial(\phi F)=\partial
F=\phi(\partial F)$. If $F\notin K$ then $\partial(\phi
F)=\partial(\sum\{(F\setminus v)\cup v_0: v\in F, (F\setminus
v)\cup v_0\in K\})$, and as we work over $\mathbb{Z}_2$, this
equals
\begin{equation}\label{eqPartialPhi}
\partial(\phi F) =
\sum\{F\setminus v : v\in F, (F\setminus v)\cup v_0\in K\} +
\end{equation}
\begin{equation}\nonumber
\sum\{(F\setminus \{u,v\})\cup v_0: u,v\in F, (F\setminus v)\cup
v_0\in K, (F\setminus u)\cup v_0\notin K\}.
\end{equation}
On the other hand $\phi(\partial F)=\phi(\sum\{F\setminus u: u\in
F, F\setminus u\in K\}) + \phi(\sum\{F\setminus u: u\in F,
F\setminus u\notin K\})$ and as we work over $\mathbb{Z}_2$, this
equals
\begin{equation}\label{eqPhiPartial}
\phi(\partial F) = \sum\{F\setminus u : u\in F, (F\setminus u)\in
K\} +
\end{equation}
\begin{equation}\nonumber
\sum\{(F\setminus \{u,v\})\cup v_0: u,v\in F, (F\setminus
\{u,v\})\cup v_0\in K, (F\setminus v)\in K,
 (F\setminus u)\notin K\}.
\end{equation}
It suffices to show that in equations (\ref{eqPartialPhi}) and
(\ref{eqPhiPartial}) the left summands on the RHSs are equal, as
well as the right summands on the RHSs. This follows from
observation \ref{obsAd} below. Thus $\phi$ is a chain map.

Second, let us check that $\phi$ is injective. Let $\pi_K$ be the
restriction map $C_*(K';\mathbb{Z}_2)\longrightarrow
\oplus\{\mathbb{Z}_2F: F\in K'\cap K\}$, $\pi_K(\sum\{\alpha_FF:
F\in K'\})=\sum\{\alpha_FF: F\in K'\cap K\}$. Similarly, let
$\pi_K^{\perp}$ be the restriction map
$C_*(K';\mathbb{Z}_2)\longrightarrow \oplus\{\mathbb{Z}_2F: F\in
K'\setminus K\}$. Note that for a chain $c\in C_*(K';\mathbb{Z}_2)$,
$c=\pi_K(c)+\pi_K^{\perp}(c)$ and $supp(\phi(\pi_K(c)))\cap
supp(\phi(\pi_K^{\perp}(c)))=\emptyset$. Assume that $c_1,c_2\in
C_*(K';\mathbb{Z}_2)$ such that $\phi(c_1)=\phi(c_2)$. Then
$\pi_K(c_1)=\phi(\pi_K(c_1))=\phi(\pi_K(c_2))=\pi_K(c_2)$, and
$\phi(\pi_K^{\perp}(c_1))=\phi(\pi_K^{\perp}(c_2))$. Note that if
$F_1,F_2 \notin K$ then $F_1,F_2\in K'$ and if $F_1\neq F_2$ then
$supp(\phi(1F_1)) \ni (F_1\setminus v_1)\cup v_0 \notin
supp(\phi(1F_2))$. Hence also
$\pi_K^{\perp}(c_1)=\pi_K^{\perp}(c_2)$. Thus $c_1=c_2$. $\square$

\begin{obs}\label{obsAd}
Let $K\mapsto K', v_0\mapsto v_1$ be an admissible contraction with
$\rm{dim}(\rm{lk}(v_0,K))\leq \rm{dim}(\rm{lk}(v_1,K))$. Let $K'\ni
F\notin K$ and $v\in F$. Then $(F\setminus v)\in K$ if and only if
$(F\setminus v)\cup v_0\in K$.
\end{obs}
$Proof$: Assume $F\setminus v\in K$. As $(F\setminus v_1)\cup v_0
\in K$ we only need to check the case $v\neq v_1$. We proceed by
induction on $\rm{dim}(F)$. As $\{v_0,v_1\}\in K$ whenever
$\rm{dim}(K)>0$ (and whenever $\rm{dim}(\rm{lk}(v_0,K))\geq 0$, if
we use the weaker local condition for admissible contractions), the
case $\rm{dim}(F)\leq 1$ is clear. (If $\rm{dim}(K)=0$ there is
nothing to prove. For the weaker local condition for admissible
contractions, if $\rm{lk}(v_0,K))=\emptyset$ then there is nothing
to prove.) By the induction hypothesis we may assume that all the
proper subsets of $(F\setminus v)\cup v_0$ are in $K$. Also
$v_0,v_1\in (F\setminus v)\cup v_0$. The admissibility of the
contraction implies that $(F\setminus v)\cup v_0\in K$. The other
direction is trivial. $\square$

\begin{lem}\label{lemInjChainMap->Smith}
Let $\phi: C_*(K';\mathbb{Z}_2)\longrightarrow C_*(K;\mathbb{Z}_2)$
be the injective chain map defined in the proof of Lemma
\ref{lemAdContr->InjChainMap} for an admissible contraction
$K\mapsto K'$. Then for every $m\geq 0$,
$\phi^*(\rm{Sm}^m([1_{K_*}]))=\rm{Sm}^m([1_{K'_*}])$ for the induced
map $\phi^*$.
\end{lem}
$Proof$: For two simplicial complexes $L$ and $L'$ and a field
$k$, the following map is an isomorphism of chain complexes:
$$\alpha=\alpha_{L,L',k}: C(L;k)\otimes_k C(L';k)\longrightarrow C(L*L';k), \ \ \alpha((1T)\otimes (1T'))=1(T\uplus T')$$
where $T\in L, T'\in L'$ and $\alpha$ is extended linearly. In
case $L=L'$ (in the definition of join we think of $L$ and $L'$ as
two disjoint copies of $L$) and $k$ is understood we denote
$\alpha_{L,L',k}=\alpha_{L}$.

Thus there is an induced chain map $\phi_*:
C_*(K'*K';\mathbb{Z}_2)\longrightarrow C_*(K*K;\mathbb{Z}_2)$,
$\phi_*=\alpha_{K}\circ \phi \otimes \phi \circ \alpha_{K'}^{-1}$
where $\phi \otimes \phi: C(K';\mathbb{Z}_2)\otimes_{\mathbb{Z}_2}
C(K';\mathbb{Z}_2)\longrightarrow
C(K;\mathbb{Z}_2)\otimes_{\mathbb{Z}_2} C(K;\mathbb{Z}_2)$ is
defined by $\phi \otimes \phi(c\otimes c')= \phi(c)\otimes \phi(c')$
(which this is a chain map).

Consider the subcomplex $C_*(K'_*;\mathbb{Z}_2)\subseteq
C_*(K'*K';\mathbb{Z}_2)$. We now verify that every $c\in
C_*(K'_*;\mathbb{Z}_2)$ satisfies $\phi_*(c)\in
C_*(K_*;\mathbb{Z}_2)$. It is enough to check this for chains of
the form $c=1(S^1\cup T^2)$ where $S,T\in K'$ and $S\cap T=
\emptyset$. For a collection of sets $A$ let $V(A)=\cup_{a\in
A}a$. Clearly if the condition
\begin{equation}\label{deljoinCond}
V(supp (\phi(S)))\cap V(supp (\phi(T)))= \emptyset
\end{equation}
is satisfied then we are done. If $v_1\notin S, v_1\notin T$, then
$\phi(S)=S, \phi(T)=T$ and (\ref{deljoinCond}) holds. If  $T\ni
v_1 \notin S$, then $\phi(S)=S$ and $V(supp \phi(T))\subseteq
T\cup\{v_0\}$. As $v_0\notin S$ condition (\ref{deljoinCond})
holds. By symmetry, (\ref{deljoinCond}) holds when $S\ni v_1
\notin T$ as well.

With abuse of notation (which we will repeat) we denote the above
chain map by $\phi$, $\phi: C_*(K'_*;\mathbb{Z}_2)\longrightarrow
C_*(K_*;\mathbb{Z}_2)$. For a simplicial complex $L$, the involution
$\tau_L:L_*\longrightarrow L_*$, $\tau_L(S^1\cup T^2)=T^1\cup S^2$
induces a $\mathbb{Z}_2$-action on $C_*(L_*;\mathbb{Z}_2)$. It is
immediate to check that $\alpha_{L,L',k}$ and $\phi \otimes \phi$
commute with these $\mathbb{Z}_2$-actions, and hence so does their
composition, $\phi$. Thus, we have proved that $\phi:
C_*(K'_*;\mathbb{Z}_2)\longrightarrow C_*(K_*;\mathbb{Z}_2)$ is a
$\mathbb{Z}_2$-chain map.

Therefore, there is an induced map on the symmetric cohomology
rings $\phi: H_S^*(K_*)\longrightarrow H_S^*(K'_*)$ which commutes
with the connecting homomorphisms $\rm{Sm}:H_S^i(L)\longrightarrow
H_S^{i+1}(L)$ for $L=K_*,K'_*$.

Let us check that for the fundamental $0$-cocycles
$\phi([1_{K_*}])=[1_{K'_*}]$ holds. A representing cochain is
$1_{K_*}: \oplus_{v\in (K_*)_0}\mathbb{Z}_2 v \longrightarrow
\mathbb{Z}_2$, $1_{K_*}(1v)=1$. As $\phi|_{C_0(K'_*)}=id$ (w.r.t.
the obvious injection $(K'_*)_0\longrightarrow (K_*)_0$), for
every $u\in (K'_*)_0$ $(\phi
1_{K_*})(u)=1_{K_*}(\phi|_{C_0(K'_*)}(u))=1_{K_*}(u)=1$, thus
$\phi(1_{K_*})=1_{K'_*}$.

As $\phi$ commutes with the Smith connecting homomorphisms, for
every $m\geq 0$, $\phi(\rm{Sm}^m(1_{K_*}))=\rm{Sm}^m(1_{K'_*})$.
$\square$

\begin{thm}\label{thmMinor->Smith}
Let $H$ and $K$ be simplicial complexes. If $H<K$ then there
exists an injective chain map
$\phi:C_*(H;\mathbb{Z}_2)\longrightarrow C_*(K;\mathbb{Z}_2)$
which induces $\phi(\rm{Sm}^m(1_{K_*}))=\rm{Sm}^m(1_{H_*})$ for
every $m\geq 0$.
\end{thm}
$Proof$: Let the sequence $K=K^0\mapsto K^1\mapsto ...\mapsto K^t=H$
demonstrate the fact that $H<K$. If $K^i\mapsto K^{i+1}$ is an
admissible contraction, then by Lemmas \ref{lemAdContr->InjChainMap}
and \ref{lemInjChainMap->Smith} it induces an injective chain map
$\phi_i: C_*(K^{i+1};\mathbb{Z}_2)\longrightarrow
C_*(K^i;\mathbb{Z}_2)$ which in turn induces
$\phi_i(\rm{Sm}^m(1_{(K^{i})_*}))=\rm{Sm}^m(1_{(K^{i+1})_*})$ for
every $m\geq 0$. If $K^i\mapsto K^{i+1}$ is a deletion - take
$\phi_i$ to be the map induced by inclusion, to obtain the same
conclusions. Thus, the composition
$\phi=\phi_0\circ...\circ\phi_{t-1}:
C_*(H;\mathbb{Z}_2)\longrightarrow C_*(K;\mathbb{Z}_2)$ is as
desired. $\square$

$Proof\ of\ Theorem\ \ref{thm o(L)}$: By Theorem
\ref{thmMinor->Smith}
$\phi(\rm{Sm}^{m}(1_{K_*}))=\rm{Sm}^{m}(1_{H_*})$. Thus if
$\rm{Sm}^{m}(1_{H_*})\neq 0$ then $\rm{Sm}^{m}(1_{K_*})\neq 0$.
 $\square$

\textbf{Remark}: The conclusion of Theorem \ref{thm o(L)} would fail
if we allow arbitrary identifications of vertices. For example, let
$K'=K_5$ and let $K$ be obtained from $K'$ by splitting a vertex
$w\in K'$ into two new vertices $u,v$, and connecting $u$ to a
non-empty proper subset of $\rm{skel}_0(K')\setminus \{w\}$, denoted
by $A$, and connecting $v$ to $(\rm{skel}_0(K')\setminus
\{w\})\setminus A$. As $K$ embeds into the $2$-sphere,
$\rm{Sm}^3(K)=0$. By identifying $u$ with $v$ we obtain $K'$, but
$\rm{Sm}^3(K')\neq 0$. To obtain from this example an example where
the edge $\{u,v\}$ is present, let $L=\rm{cone}(K)\cup \{u,v\}$
($\rm{cone}(K)$ is the cone over $K$), and let $L'$ be the complex
obtained form $L$ by identifying $u$ with $v$. Then $\rm{Sm}^4(L)=0$
while $\rm{Sm}^4(L')\neq 0$.

\begin{ex}\label{exFindMinor} Let $K$ be the simplicial complex
spanned by the following collection of $2$-simplices:
$(\binom{[7]}{3}\setminus \{127,137,237 \})\cup
\{128,138,238,178,278,378 \}$.
\end{ex}
$K$ is not a subdivision of $H(3)$, and its geometric realization
even does not contain a subspace homeomorphic to $H(3)$ (as there
are no $7$ points in $|K|$, each with a neighborhood whose
boundary contains a subspace which is homeomorphic to $K_6$).
Nevertheless, contraction of the edge $78$ is admissible and
results in $H(3)$. By Theorem \ref{thm o(L)} $K$ has a
non-vanishing Van Kampen's obstruction in dimension $5$, and hence
is not embeddable in the $4$-sphere.

\begin{ex}\label{exUli}
Let $K_1$ be a triangulation of $S^1$ (the 1-sphere) and let $K_2$
be a triangulation of $S^2$. Then $K=K_1*K_2$ is a triangulation
of $S^4$. Let $T$ be a missing triangle of $K$ and
$L=\rm{skel}_2(K)\cup\{T\}$. Then $L$ does not embed in
$\mathbb{R}^4$.
\end{ex}
$Proof$: It is known and easy to prove that every $2$-sphere may be
reduced to the boundary of the tetrahedron by a sequence of
admissible contractions in a way that fixes a chosen triangle from
the original triangulation (e.g. \cite{Whiteley}, Lemma 6). This
guarantees the existence of sequences of admissible contractions as
described below.

Case 1: $\partial(T)=K_1$. There exists a sequences of admissible
contractions (of vertices from $K_2$) which reduces $L$ to $H(3)$.
By Theorems \ref{thm o(L)}, \ref{thmSmithH(d)} and
\ref{thmSmith->nonEmb}, $L$ does not embed in $\mathbb{R}^4$.

Case 2: $\partial(T)\neq K_1$. Hence $\partial(T)\subseteq K_2$
and separates $K_2$ into two disks. By performing admissible
contractions of pairs of vertices within each of these disks, and
within $K_1$, we can reduce $L$ to the $2$-skeleton of the join
$L_1*L_2$ where $L_1$ is the boundary of a triangle and $L_2$ is
two boundaries of tetrahedra glued along a triangle. Let $v$ be a
vertex which belongs to exactly one of the two tetrahedra which
were used to define $L_2$. Deleting $v$ from $L$ results in $H(3)$
minus one triangle which consists of the vertices of $L_1$. Hence
the subcomplex $L'=(L-v)\cup (L_1*\{v\})$ of $L$ is admissibly
contracted into $H(3)$ by contracting an edge which contains $v$.
Thus, $H(3)<L$ and by Theorems \ref{thm o(L)}, \ref{thmSmithH(d)}
and \ref{thmSmith->nonEmb}, $L$ does not embed in $\mathbb{R}^4$.
$\square$

Example \ref{exUli} is a special case of the following conjecture,
a work in progress of Uli Wagner and the author.

\begin{conj}\label{conjUli}
Let $K$ be a triangulated $2d$-sphere and let $T$ be a missing
$d$-face in $K$. Let $L=\rm{skel}_d(K)\cup\{T\}$. Then $L$ does
not embed in $\mathbb{R}^{2d}$.
\end{conj}

\section{The obstruction over $\mathbb{Z}$}\label{secOverZ}
More commonly in the literature, Van Kampen's obstruction is defined
via deleted products and with $\mathbb{Z}$ coefficients, where,
except for $2$-simplicial complexes, its vanishing is also
sufficient for embedding of the complex in a Euclidean space of
double its dimension. We obtain an analogue of Theorem \ref{thm
o(L)} for this context.

The presentation of the background on the obstruction here is based
on the ones in \cite{Novik}, \cite{Wu} and \cite{Ummel}.

Let $K$ be a finite simplicial complex. Its deleted product is
$K\times K \setminus \{(x,x): x\in K\}$, employed with a fixed-point
free $\mathbb{Z}_2$-action $\tau(x,y)=(y,x)$. It
$\mathbb{Z}_2$-deformation retracts into $K_{\times}=\cup\{S\times
T: S,T\in K, S\cap T=\emptyset\}$, with which we associate a cell
chain complex over $\mathbb{Z}$: $C_{\bullet}(K_{\times})=\bigoplus
\{\mathbb{Z}(S\times T): S\times T \in K_\times\}$ with a boundary
map $\partial(S\times T)=\partial S \times T + (-1)^{\rm{dim} S}
S\times
\partial T$, where $S\times T$ is a $\rm{dim}(S\times T)$-chain.
The dual cochain complex consists of the $j$-cochains
$C^j(K_{\times})=\rm{Hom}_{\mathbb{Z}}(C_j(K_{\times}),\mathbb{Z})$
for every $j$.

There is a $\mathbb{Z}_2$-action on $C_{\bullet}(K_{\times})$
defined by $\tau(S\times T)=(-1)^{\rm{dim}(S)\rm{dim}(T)}T\times S$.
As it commutes with the coboundary map, by restriction of the
coboundary map we obtain the subcomplexes of symmetric cochains
$C_s^{\bullet}(K_{\times})=\{c\in C^{\bullet}(K_{\times}):
\tau(c)=c\}$ and of antisymmetric cochains
$C_a^{\bullet}(K_{\times})=\{c\in C^{\bullet}(K_{\times}):
\tau(c)=-c\}$. Their cohomology rings are denoted by
$H_s^{\bullet}(K_{\times})$ and $H_a^{\bullet}(K_{\times})$
respectively. Let $H_{\rm{eq}}^{m}$ be $H_s^{m}$ for $m$ even and
$H_a^{m}$ for $m$ odd.

For every finite simplicial complex $K$ there is a unique
$\mathbb{Z}_2$-map, up to $\mathbb{Z}_2$-homotopy, into the infinite
dimensional sphere $i: K_{\times}\rightarrow S^{\infty}$, and hence
a uniquely defined map $i^*:
H_{\rm{eq}}^{\bullet}(S^{\infty})\rightarrow
H_{\rm{eq}}^{\bullet}(K_{\times})$. For $z$ a generator of
$H_{\rm{eq}}^{m}(S^{\infty})$ call
$o^m=o^m_{\mathbb{Z}}(K_{\times})=i^*(z)$ the Van Kampen
obstruction; it is uniquely defined up to a sign. It turns out to
have the following explicit description: fix a total order $<$ on
the vertices of $K$. It evaluates elementary symmetric chains of
even dimension $2m$ by
\begin{equation}\label{eqVKs}
o^{2m}((1+\tau)(S\times T))= \{^{1\ \rm{if}\ \rm{the}\
\rm{unordered}\ \rm{pair}\ \{S,T\}\ \rm{is}\ \rm{of}\ \rm{the}\
\rm{form}\ s_0<t_0<..<s_m<t_m} _{0 \ \rm{for}\ \rm{other}\
\rm{pairs}\ \{S,T\}}
\end{equation}
and evaluates elementary antisymmetric chains of odd dimension
$2m+1$ by
\begin{equation}\label{eqVKa}
o^{2m+1}((1-\tau)(S\times T))= \{^{1\ \rm{if}\ \{S,T\}\ \rm{is}\
\rm{of}\ \rm{the}\ \rm{form}\ t_0<s_0<t_1<..<t_m<s_m<t_{m+1}} _{0 \
\rm{for}\ \rm{other}\ \rm{pairs}\ \{S,T\}}
\end{equation}
where the $s_l$'s are elements of $S$ and the $t_l$'s are elements
of $T$. Its importance to embeddability is given in the following
classical result:
\begin{thm}\label{thmVK}\cite{VanKampen,Shapiro,Wu}
If a simplicial complex $K$ embeds in $\mathbb{R}^m$ then
$H_{\rm{eq}}^{\bullet}(K_{\times})\ni
o^m_{\mathbb{Z}}(K_{\times})=0$. If $K$ is $m$-dimensional and
$m\neq 2$ then $o^{2m}_{\mathbb{Z}}(K_{\times})=0$ implies that $K$
embeds in $\mathbb{R}^{2m}$.
\end{thm}

In relation to higher minors, the analogue of Theorem \ref{thm o(L)}
holds:
\begin{thm}\label{thm o_Z}
Let $H$ and $K$ be simplicial complexes. If $H<K$ and
$o^m_{\mathbb{Z}}(H_{\times})\neq 0$ then
$o^m_{\mathbb{Z}}(K_{\times})\neq 0$.
\end{thm}
From Theorems \ref{thm o_Z} and \ref{thmVK} it follows that
Conjecture \ref{conjMinorEmb} is true when $2\rm{dim}(H) = m \neq 4$
(and, trivially, when $2\rm{dim}(H) < m$).

$Proof\ of\ Theorem\ \ref{thm o_Z}$: Fix a total order on the
vertices of $K$, $v_0<v_1<..<v_n$ and consider an admissible
contraction $K\mapsto K'$ where $K'$ is obtained from $K$ by
identifying $v_0\mapsto v_1$ (shortly this will be shown to be
without loss of generality). Define a map $\phi$ as follows: for
$F\in K'$
\begin{equation}\label{phiTimes}
\phi(F)=\{^{F \ \ \rm{if}\ F\in K}_ {\sum\{\rm{sgn}(v,F)(F\setminus
v)\cup v_0:\ v\in F, (F\setminus v)\cup v_0\in K\}\ \ \rm{if}\
F\notin K}
\end{equation}
where $\rm{sgn}(v,F)=(-1)^{|\{t\in F: t<v\}|}$. Extend linearly to
obtain an injective $\mathbb{Z}$-chain map $\phi:
C_{\bullet}(K')\longrightarrow C_{\bullet}(K)$. (The check that this
map is indeed an injective $\mathbb{Z}$-chain map is similar to the
proof of Lemma \ref{lemAdContr->InjChainMap}.) In case we contract a
general $a\mapsto b$, for the signs to work out consider the map
$\tilde{\phi}=\pi^{-1}\phi\pi$ rather than $\phi$, where $\pi$ is
induced by a permutation on the vertices which maps $\pi(a)=v_0,\
\pi(b)=v_1$. Then $\tilde{\phi}$ is an injective $\mathbb{Z}$-chain
map.

As $\phi(S\times T):=\phi(S)\times \phi(T)$ commutes with the
$\mathbb{Z}_2$ action and with the boundary map on the chain complex
of the deleted product, $\phi$ induces a map
$H_{\rm{eq}}^{\bullet}(K_{\times})\rightarrow
H_{\rm{eq}}^{\bullet}(K'_{\times})$. It satisfies
$\phi^*(o^m_{\mathbb{Z}}(K_{\times}))=o^m_{\mathbb{Z}}(K'_{\times})$
for all $m\geq 1$. The checks are straightforward (for proving the
last statement, choose a total order with contraction which
identifies the minimal two elements $v_0\mapsto v_1$, and show
equality on the level of cochains). We omit the details.

If $K\mapsto K'$ is a deletion, consider the injection $\phi:
K'\rightarrow K$ to obtain again an induced map with
$\phi^*(o^m_{\mathbb{Z}}(K_{\times}))=o^m_{\mathbb{Z}}(K'_{\times})$.

Let the sequence $K=K^0\mapsto K^1\mapsto ...\mapsto K^t=H$
demonstrate the fact that $H<K$. By composing the corresponding maps
as above we obtain a map $\phi^*$ with
$\phi^*(o^m_{\mathbb{Z}}(K_{\times}))=o^m_{\mathbb{Z}}(H_{\times})$
and the result follows. $\square$

\section{Topology preserving edge contractions}\label{secTop}
\subsection{PL manifolds}
$Proof\ of\ Theorem\ \ref{thmDey}$: Let $M$ be a PL-triangulation
of a compact $d$-manifold without boundary. Let $ab$ be an edge of
$M$ and let $M'$ be obtained from $M$ by contracting $a\mapsto b$.
We will prove that if the Link Condition (\ref{eqLinkCond}) holds
for $ab$ then $M$ and $M'$ are PL-homeomorphic, and otherwise they
are not homeomorphic (not even 'locally homologic'). For $d=1$ the
assertion is clear. Assume $d>1$. Denote the $closed$ star of a
vertex $a$ in $M$ by $\rm{st}(a,M)= \{T\in M: T\cup \{a\}\in M\}$
and denote its antistar by $\rm{ast}(a,M)=  \{T\in M: a\notin
T\}$.

Denote $B(b)=\{b\}*\rm{ast}(b,\rm{lk}(a,M))$ and
$L=\rm{ast}(a,M)\cap B(b)$. Then $M'=\rm{ast}(a,M)\cup_{L}B(b)$.
As $M$ is a PL-manifold without boundary, $\rm{lk}(a,M)$ is a
$(d-1)$-PL-sphere (see e.g. Corollary 1.16 in \cite{Hudson}). By
Newman's theorem (e.g. \cite{Hudson}, Theorem 1.26)
$\rm{ast}(b,\rm{lk}(a,M))$ is a $(d-1)$-PL-ball. Thus $B(b)$ is a
$d$-PL-ball. Observe that $\partial
(B(b))=\rm{ast}(b,\rm{lk}(a,M))\cup \{b\} *
\rm{lk}(b,\rm{lk}(a,M)) = \rm{lk}(a,M) = \partial (\rm{st}(a,M))$.

The identity map on $\rm{lk}(a,M)$ is a PL-homeomorphism $h:
\partial (B(b))\rightarrow \partial (\rm{st}(a,M))$, hence it
extends to a PL-homeomorphism $\tilde{h}: B(b)\rightarrow
\rm{st}(a,M)$ (see e.g. \cite{Hudson}, Lemma 1.21).

Note that $L=\rm{lk}(a,M)\cup (\{b\}*(\rm{lk}(a,M)\cap
\rm{lk}(b,M)))$.

If $\rm{lk}(a)\cap \rm{lk}(b)=\rm{lk}(ab)$ (in $M$) then
$L=\rm{lk}(a,M)$, hence gluing together the maps $\tilde{h}$ and
the identity map on $\rm{ast}(a,M)$ results in a PL-homeomorphism
from $M'$ to $M$.

If $\rm{lk}(a)\cap \rm{lk}(b)\neq \rm{lk}(ab)$ (in $M$) then
$\rm{lk}(a,M)\subsetneqq L$. The case $L=B(b)$ implies that
$M'=\rm{ast}(a,M)$ and hence $M'$ has a nonempty boundary, showing
it is not homeomorphic to $M$. A small punctured neighborhood of a
point in the boundary of $M'$ has trivial homology while all small
punctured neighborhoods of points in $M$ has non vanishing
$(d-1)$-th homology. This is what we mean by 'not even locally
homologic': $M$ and $M'$ have homologically different sets of
small punctured neighborhoods.

We are left to deal with the case $\rm{lk}(a,M)\subsetneqq
L\subsetneqq B(b)$. As $L$ is closed there exists a point $t\in
L\cap \rm{int}(B(b))$ with a small punctured neighborhood $N(t,M')$
which is not contained in $L$. For a subspace $K$ of $M'$ denote by
$N(t,K)$ the neighborhood in $K$ $N(t,M')\cap K$. Thus
$N(t,M')=N(t,\rm{ast}(a,M))\cup_{N(t,L)}N(t,B(b))$. We get a
Mayer-Vietoris exact sequence in reduced homology:
\begin{equation}\label{eqMV}
H_{d-1}N(t,L)\rightarrow H_{d-1}N(t,\rm{ast}(a,M))\oplus
H_{d-1}N(t,B(b))\rightarrow H_{d-1}N(t,M')\rightarrow
\end{equation}
\begin{equation}\nonumber
H_{d-2}N(t,L)\rightarrow H_{d-2}N(t,\rm{ast}(a,M))\oplus
H_{d-2}N(t,B(b)).
\end{equation}
Note that $N(t,\rm{ast}(a,M))$ and $N(t,B(b))$ are homotopic to
their boundaries which are $(d-1)$-spheres. Note further that
$N(t,L)$ is homotopic to a proper subset $X$ of $\partial
(N(t,B(b)))$ such that the pair $(\partial (N(t,B(b))),X)$ is
triangulated. By Alexander duality $H_{d-1}N(t,L)=0$. Thus,
(\ref{eqMV}) simplifies to the exact sequence
$$0\rightarrow \mathbb{Z}\oplus \mathbb{Z}\rightarrow H_{d-1}N(t,M')\rightarrow
H_{d-2}N(t,L)\rightarrow 0.$$ Thus, $\rm{rank}(H_{d-1}N(t,M'))\geq
2$, hence $M$ and $M'$ are not locally homologic, and in
particular are not homeomorphic. $\square$

\textbf{Remarks}: (1) Omitting the assumption in Theorem
\ref{thmDey} that the boundary is empty makes both implications
incorrect. Contracting an edge to a point shows that the Link
Condition is not sufficient. Contracting an edge on the boundary
of a cone over an empty triangle shows that the Link Condition is
not necessary.

(2) The necessity of the Link Condition holds also in the
topological category (and not only in the PL category), as the
proof of Theorem \ref{thmDey} shows. Indeed, for this part we only
used the fact that $B(b)$ is a pseudo manifold with boundary
$\rm{lk}(a,M)$ (not that it is a ball); taking the point $t$ to
belong to exactly two facets of $B(b)$. The following part, in the
topological category, is still open:

\begin{prob}\label{probTop}
Given an edge in a triangulation of a compact manifold without
boundary which satisfies the Link Condition, is it true that its
contraction results in a homeomorphic space? Or at least in a
space of the same homotopic or homological type?
\end{prob}
A Mayer-Vietoris argument shows that such topological manifolds $M$
and $M'$ have the same Betti numbers; both $\rm{st}(a,M)$ and $B(b)$
are cones and hence their reduced homology vanishes.

A candidate for a counterexample for Problem \ref{probTop} may be
the join $M=T*P$ where $T$ is the boundary of a triangle and $P$ a
triangulation of Poincar\'{e} homology $3$-sphere, where an edge
with one vertex in $T$ and the other in $P$ satisfies the Link
Condition. By the double-suspension theorem (Edwards
\cite{Edwards} and Cannon \cite{Cannon}) $M$ is a topological
$5$-sphere.

Walkup \cite{Walkup} mentioned, without details, the necessity of
the Link Condition for contractions in topological manifolds, as
well as the sufficiency of the Link Condition for the $3$
dimensional case (where the category of PL-manifolds coincides
with the topological one); see \cite{Walkup}, p.82-83.

\subsection{PL spheres}\label{subsecSpheres}
In this section we use some terminology from $f$-vectors theory;
readers unfamiliar with this terminology can consult
\cite{StanleyBook}.

\begin{defin}\label{defStrongContr}
Boundary complexes of simplices are \emph{strongly edge
decomposable} and, recursively, a triangulated PL-manifold $S$ is
\emph{strongly edge decomposable} if it has an edge which satisfies
the Link Condition (\ref{eqLinkCond}) such that both its link and
its contraction are strongly edge decomposable.
\end{defin}
By Theorem \ref{thmDey} the complexes in Definition
\ref{defStrongContr} are all triangulated PL-spheres. Note that
every $2$-sphere is strongly edge decomposable.

Let $vu$ be an edge in a simplicial complex $K$ which satisfies the
Link Condition, whose contraction $u\mapsto v$ results in the
simplicial complex $K'$. Note that the $f$-polynomials satisfy
$$f(K,t)=f(K',t)+t(1+t)f(\rm{lk}(\{vu\},K),t),$$
hence the $h$-polynomials satisfy
\begin{equation}\label{eq-hContraction}
h(K,t)=h(K',t)+t h(\rm{lk}(\{vu\},K),t).
\end{equation}
We conclude the following:
\begin{cor}
The $g$-vector of strongly edge decomposable triangulated spheres is
non negative. $\square$
\end{cor}
Is it also an $M$-sequence? The strongly edge decomposable spheres
(strictly) include the family of triangulated spheres which can be
obtained from the boundary of a simplex by repeated Stellar
subdivisions (at any face); the later are polytopal, hence their
$g$-vector is an $M$-sequence. For the case of subdividing only at
edges (\ref{eq-hContraction}) was considered by Gal (\cite{Gal},
Proposition 2.4.3).

\section{Graph minors versus higher minors}\label{secVersus}
While Theorem \ref{thm o(L)} is an instance of a property of graph
minors which generalizes to higher minors, this is not always the
case. Let us mention some properties which do not generalize, and
others for which we do not know whether they generalize or not.

\begin{itemize}
\item For graphs, if $K$ is a subdivision of $H$ then $H$ is a
minor of $K$. This is not the case for higher minors.

\begin{ex}
Let $H$ be a triangulated PL $3$-sphere whose triangulation contains
a knotted triangle $\{12,23,13\}$ (e.g. \cite{Lutz2} for an example
with few vertices and references to Hachimori's first examples. In
\cite{Hachimori-Ziegler} such spheres were proved to be
non-constructible). Then $H$ is a subdivision of $\partial\Delta^4$,
the boundary complex of the $4$-simplex, but $\partial\Delta^4$ is
not a minor of $H$.
\end{ex}
$Proof$: Consider, by contradiction, a sequence of deletions and
admissible contractions starting at $H$ and ending at
$\partial\Delta^4$. Any deletion would result in a complex with a
vanishing $3$-homology; further deletions and contractions would
keep the $3$-homology being zero as they induce the injective chain
map from Theorem \ref{thmMinor->Smith}. Thus the sequence contains
only contractions. Any admissible contraction, assuming we haven't
reached $\partial\Delta^4$ yet, must satisfy the Link Condition
(\ref{eqLinkCond}) - as by Alexander duality a sphere can not
contain a sphere of the same dimension as a proper subspace. If a
contraction $a\mapsto b$ satisfies $a\neq 1,2,3$, the
PL-homeomorphism constructed in the proof of Theorem \ref{thmDey}
shows that it results in a PL $3$-sphere with $\{12,23,13\}$ a
knotted triangle. It suffices to show that a contraction where $a\in
\{1,2,3\}$ also results in a triangulation with a knotted triangle,
as this would imply that $\partial\Delta^4$ can never be reached, a
contradiction. Without loss of generality $a=1$. As $\{12,23,13\}$
is knotted in $M$, the Link Condition implies $b\neq 2,3$ and
$\{b,2,3\}\notin M$. Thus $\{b2,23,b3\}$ is an induced subcomplex in
$M'$, and hence there is a deformation retract of $M'-\{b2,23,b3\}$
onto the induced subcomplex
$M'[V(M')-\{b,2,3\}]=M[V(M)-\{b,1,2,3\}]$, where $V(K)$ is the set
of vertices of a complex $K$ (e.g. \cite{Dancis}, Lemma 4).
Similarly, $M[V(M)-\{b,1,2,3\}]$ is a deformation retract of
$M[V(M)-\{1,2,3\}]-\{b\}$. To show that the fundamental group
$\pi_1(M'-\{b2,23,b3\})\neq 0$ we will show that
$\pi_1(M[V(M)-\{1,2,3\}]-\{b\})\neq 0$. We use Van Kampen's theorem
for the union $M[V(M)-\{1,2,3\}]=(M[V(M)-\{1,2,3\}]-\{b\})\cup
int(star(b,M[V(M)-\{1,2,3\}]))$: note that the intersection is a
deformation retract of $lk(b,M)$ minus the induced subcomplex on
$\{1,2,3\}$ in it, which is path-connected and simply connected. We
conclude that $\pi_1(M[V(M)-\{1,2,3\}]-\{b\})\cong
\pi_1(M[V(M)-\{1,2,3\}]) \neq 0$, as $\{12,23,13\}$ is knotted in
$M$. $\square$

\item For a graph $K$ on $n$ vertices, if $K$ has more than $3n-6$
edges then it contains a $K_5$ minor (Mader proved that it even
contains a $K_5$ subdivision \cite{Mader2}). Is the following
generalization to higher minors true?:

\begin{prob}\label{probC(d,n)}
Let $C(d,n))$ be the boundary complex of a cyclic $d$-polytope on
$n$ vertices, and let $K$ be a simplicial complex on $n$ vertices.
Does $f_{d}(K)>f_{d}(C(2d+1,n))$ imply $H(d+1)<K$ ?
\end{prob}
\begin{ex}
Let $M_L$ be the vertex transitive neighborly $4$-sphere on $15$
vertices $\rm{manifold}\_(4,15,5,1)$ found by Frank Lutz
\cite{Lutz}.

$M_L$ has no universal edges, i.e., every edge is contained in a
missing triangle.
\end{ex}

It is possible that $K$ equals the $2$-skeleton of $M_L$ union with
a missing triangle would provide a counterexample to Problem
\ref{probC(d,n)}.

\item If $K$ is the graph of a triangulated $2$-sphere union with
a missing edge then it contains a $K_5$ minor (the condition implies
having more than $3n-6$ edges). Is the following generalization to
higher minors true?:

\begin{prob}
Let $K$ be the union of the $d$-skeleton of a triangulated
$2d$-sphere with a missing $d$-face. Does $H(d+1)<K$ ?
\end{prob}
It is possible that $K$ equals the onion of the $2$-skeleton of
$M_L$ with a missing triangle would provide a counterexample. But if
true, then by Theorems \ref{thm o(L)}, \ref{thmSmithH(d)} and
\ref{thmSmith->nonEmb}, Conjecture \ref{conjUli} will follow.

\item A Robertson-Seymour type theorem does not hold for embeddability in higher dimensional
spheres:
\begin{prop}
For any $d\geq 2$ There exist infinitely many $d$-complexes not
embeddable in the $2d$-sphere such all of their proper minors do
embed in the $2d$-sphere.
\end{prop}
$Proof$: By identifying disjoint pair of points, each pair to a
point, where each pair lies in the interior of a facet of $H(d+1)$,
one obtains topological spaces which are not embeddable in the
$2d$-sphere but such that any proper subspace of them is. This was
proved by Zaks \cite{Zaks} for $d>2$ and later by Ummel \cite{Ummel}
for $d=2$. By choosing say $m$ such pairs in each facet, one obtains
infinitely many pairwise non-homeomorphic such spaces when $m$
varies. To conclude the claim it suffices to triangulate these
spaces in a way that no contraction would be admissible; this is
indeed possible
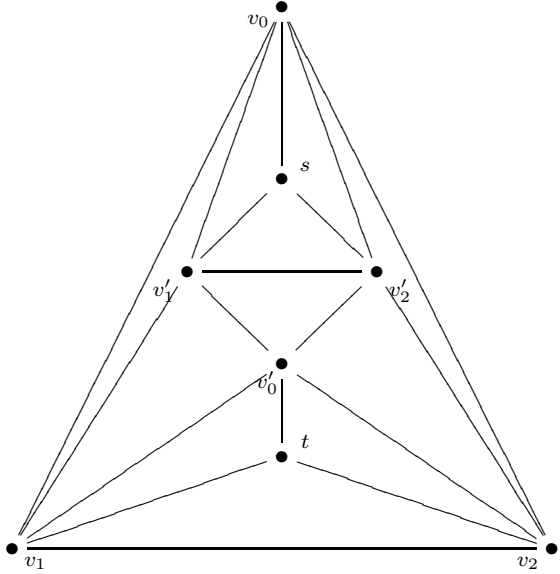
\begin{figure}\label{FigSd}
\newcommand{\edge}[1]{\ar@{-}[#1]}
\newcommand{\uulab}[1]{\ar@{}[]^<<{#1}}
\newcommand{\lulab}[1]{\ar@{}[l]^<<{#1}}
\newcommand{\rulab}[1]{\ar@{}[r]^<<{#1}}
\newcommand{\ldlab}[1]{\ar@{}[l]^<<{#1}}
\newcommand{\rdlab}[1]{\ar@{}[r]_<<{#1}}
\newcommand{\node}{*+[O][F-]{ }}
\centerline{ \xymatrix{
 & & & \bullet \lulab{v_0} \edge{ddddddlll} \edge{dddl} \edge{dd} \edge{dddr} \edge{ddddddrrr} & & & \\
& & & & & & & \\
& & & \bullet \rulab{s} \edge{dr} \edge{dl} & & &\\
& & \bullet \lulab{v'_1} \edge{rr} \edge{dddll} & & \bullet \rdlab{v'_2} \edge{dddrr} & & \\
& & & \bullet \uulab{v'_0} \edge{ul} \edge{ur} \edge{ddlll} \edge{ddrrr} & & & \\
& & & \bullet \rulab{t} \edge{u} \edge{dlll} \edge{drrr} & & & \\
 \bullet \rdlab{v_1} \edge{rrrrrr} & & & & & &\bullet \lulab{v_2}
 }} \caption {Subdivision of a small facet $F=\{v_0,v_1,v_2\}$.}
\end{figure}
(see Figure 1 for an illustration): first subdivide each facet into
$m$ small facets say. To identify simplicialy a pair of points $s,t$
in the interior of a small facet $F=\{v_0,...,v_d\}$ first further
subdivide $F$ as follows. Consider the prism $[0,1]\times
\{v_1,..,v_d\}$ with bottom $\{v_1,..,v_d\}$ and top
$\{v'_1,..,v'_d\}$ and triangulate the cylinder $[0,1]\times
\partial\{v_1,..,v_d\}$ without adding new vertices (this is standard).
Now cone with a
vertex $v'_0$ over $\partial([0,1]\times \{v_1,..,v_d\})$ to obtain
a triangulation of the prism, and further cone with the vertex $v_0$
over $\partial([0,1]\times \{v_1,..,v_d\})-\{v_1,..,v_d\}$ to
obtain, together with the prism, a triangulation of $F$. Subdivide
$\{v_1,..,v_d,v'_0\}$ by staring from a vertex $s$ in its interior,
and subdivide $\{v'_1,..,v'_d,v_0\}$ by staring from a vertex $t$ in
its interior. Note that identifying $s\mapsto t$ results in a
complex where each pair of vertices from $v_0,..,v_d,v'_0,..,v'_d,t$
is contained in a missing face of dimension $<d$ (a facet for a pair
from $v_0,..,v_d$ or from $v'_0,..,v'_d$, and an edge or a triangle
with the vertex $t$ for the rest of the pairs). $\square$

\end{itemize}

\section*{Acknowledgments}
I wish to thank Gil Kalai and Eric Babson for helpful discussions,
and Uli Wagner and Karanbir Sarkaria for helpful remarks on the
presentation. Part of this work was done during the author's stay at
Institut Mittag-Leffler, supported by the ACE network.

\end{document}